\documentclass[a4j,12pt]{article}
\usepackage{graphicx,color}
\usepackage{amsthm}
\usepackage{amsmath}
\usepackage{amssymb}
\usepackage{adjustbox}

\newtheorem{theorem}{Theorem}

\begin{document}

\begin{center}

{\bf \Large ChatGPT solved the dynamic construction problem of Malfatti circles}

\end{center}

\begin{center}
{\large  Kazushi Ahara} \\
Meiji University\\
ahara@meiji.ac.jp\\

\bigskip

\begin{minipage}{13cm}
{\small {\bf abstract}: In this paper, we address the dynamic construction problem for Malfatti circles with respect to a fixed triangle, using the LC and CC modules. We show that the configuration of the Malfatti circles constitutes a stable fixed point and provide lower bounds for the main parameters of the LC and CC modules. We now remark that all mathematical proofs in this paper were generated by ChatGPT.}
\end{minipage}
\end{center}

\bigskip

{\small keyword: ChatGPT, dynamic construction, Malfatti circles}

\section{Introduction}

\subsection{Before introduction}
This paper presents a problem-solving exercise that ChatGPT performed on a mathematical topic conceived by the author. The contributions of ChatGPT and the author to the paper are as follows: ChatGPT was responsible for constructing the proof and drafting the text. The author’s responsibilities included organizing the chapters, deciding on the title, providing prompts, revising the draft, and checking the wording. We leave it to the reader to decide whether this form of authorship is “normal.”

About related researches, see Section \ref{sec:related}

\subsection{Dynamic construction of the planer geometry}
The author developed a plane geometry construction software package called PointLine \cite{ahara1}.
In this software, all construction modules are generated by a Markov process,
and the mathematical goal is to show that these modules converge to the desired figure by repeatedly applying the Markov process.
We will refer to the process of obtaining a figure as the limit of the convergence as {\bf dynamic construction}.

The types of construction modules prepared in PointLine include \lq midpoint\rq, \lq point on a line\rq, \lq circle tangent to a line\rq, and \lq circles tangent to each other\rq.
While its user interface is comparable to that of interactive geometry software such as GeoGebra,
the underlying principles that make the figures valid are entirely different.

Here, we consider the following problems.
Let ABC be a given non-degenerate triangle.
We construct three distinct circles in the same plane,
using only the modules \lq a line is tangent to a circle' and \lq a circle is tangent to a circle\rq,
to form a set of modules
such that the configuration of the Malfatti circles of the given triangle is a fixed point.
In this case, is the Malfatti circle a stable fixed point?
Furthermore, what is the permissible range of the module parameters?
And can the global initial-value problem be solved?

In this paper, we report that ChatGPT has completely solved the first problem,
has usefully resolved the second problem,
and has presented substantive insights regarding the third problem.

\subsection{Malfatti circles}

We will explain Malfatti circles, which served as the starting point for this problem.
For a given triangle $\mathrm{ABC}$, three circles $C_1$, $C_2$, and $C_3$ are said to be Malfatti circles if they satisfy the following conditions: $C_1$ and $C_2$ are tangent to line segment $\mathrm{AB}$, $C_2$ and $C_3$ are tangent to line segment $\mathrm{BC}$, $C_3$ and $C_1$ are tangent to line segment $\mathrm{CA}$, and $C_1$, $C_2$, and $C_3$ are mutually circumscribed.
In general, the method for constructing Malfatti circles using a compass and ruler, as described by Jakob Steiner, is well known.
\begin{figure}[hbt]
\begin{center}
\includegraphics[height=4cm]{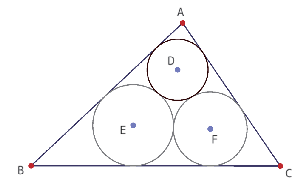}\\
\caption{Malfatti circles }
\label{fig:ajima}
\end{center}
\end{figure}

\subsection{LC (Line–Circle) Module of PointLine}

In this subsection, we will introduce the definition of the LC module in PointLine. 
Let $L$ be a line and $C$ a circle with center $\mathrm{P}$ and radius $r$. In this paper, we assume that $L$ is fixed.
We will define the LC-module of a parameter $ q$ as follows.
Let
$h=\operatorname{dist}(\mathrm{P}, L)$
be the distance from the center of the circle to the line. Define the signed contact error by
$d=r-h.$
Thus, $d=0$ exactly when the circle is tangent to the line.

The purpose of the LC module is to reduce this contact error while allowing both the center and the radius of the circle to change.
Let $\mathbf n$ be the unit vector perpendicular to $L$, pointing from $L$ toward $\mathrm{P}$. In one application of the LC module, the correction is divided equally between the motion of the center and the change of the radius:
\[
\mathrm{P}' = \mathrm{P}+\dfrac{(1-{q})d}{2}\mathbf{n}, \qquad r'=r-\frac{(1-{q})d}{2}.
\]
Since the center moves by ${(1-{q})d}/{2}$ in the normal direction, its distance from the line changes from $h$ to $h'=h+{(1-{q})d}/{2}$.
Therefore the new contact error is
\[
d'
=r'-h'
=(r-h)-(1-{q})d
={q} d.
\]
Hence one application of the LC module multiplies the contact error by ${q}$.
If $d>0$, the radius is too large relative to the distance from the line. The LC module therefore decreases the radius and moves the center away from the line.
If $d<0$, the radius is too small relative to the distance from the line. The same formulas automatically increase the radius and move the center toward the line.
Thus, no separate treatment of the two cases is necessary.
If $d_0$ is the initial contact error, then after $n$ applications of the module, $d_n=q^nd_0$. Thus, for $|q|<1$, we have $d_n\to 0$ as $n\to\infty$.
The line and the circles converge asymptotically to a tangent configuration. 

\begin{figure}[hbt]
\begin{center}
\includegraphics[height=5cm]{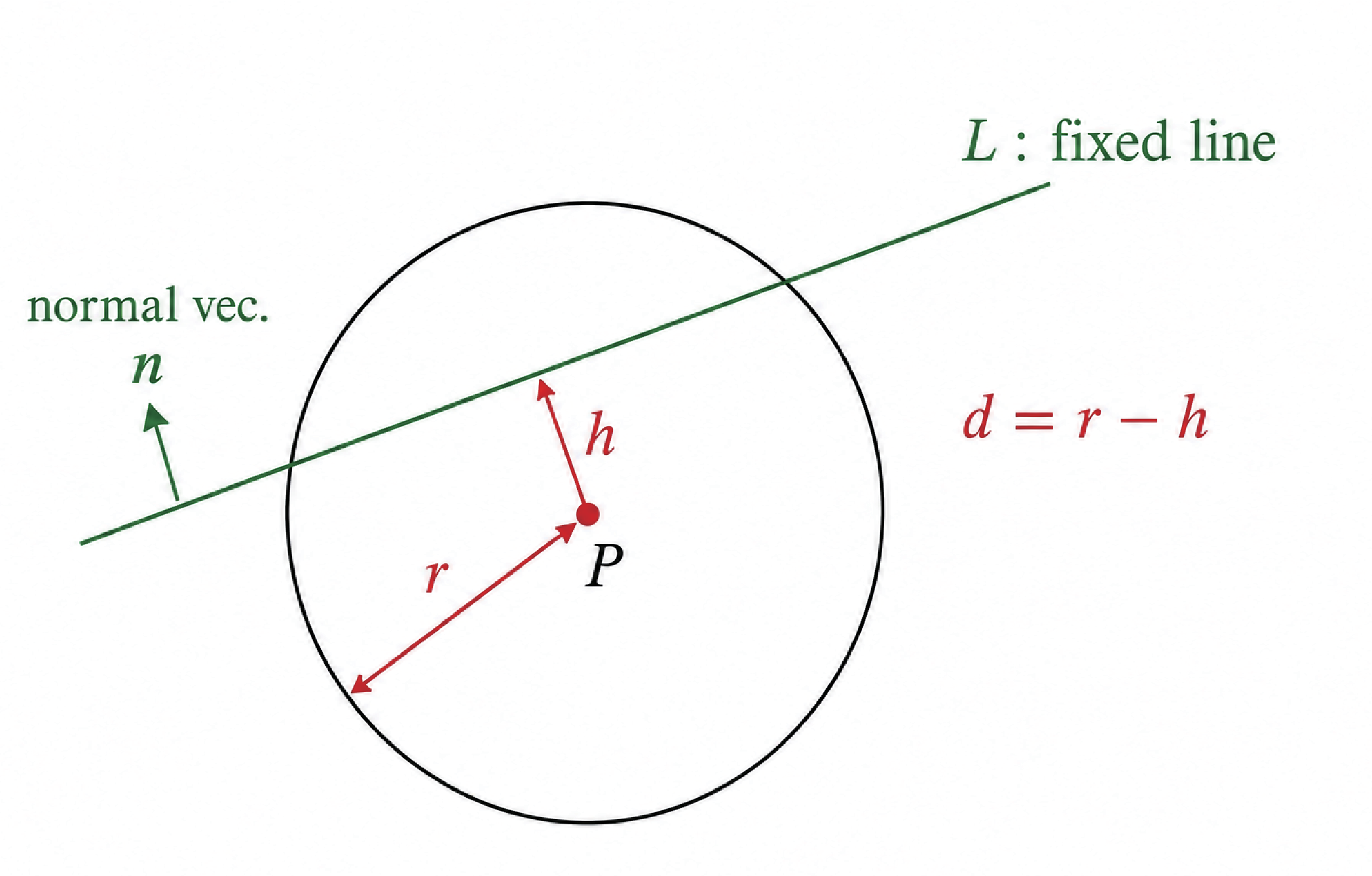}
\caption{LC Module}
\label{fig:LCmodule}
\end{center}
\end{figure}

\subsection{CC (Circle–Circle-circumscribed) Module of PointLine}

In this subsection, we will introduce the definition of the CC module in PointLine. 
Let $C_1$ and $C_2$ be two circles with centers $\mathrm{P}_1$ and $\mathrm{P}_2$, and radii $r_1$ and $r_2$, respectively. 
Let $R=\|\mathrm{P}_2-\mathrm{P}_1\|$ be the distance between the two centers, and suppose that $R\neq 0$. 
Define the signed contact error by $d=R-r_1-r_2$.
Thus, $d=0$ exactly when the two circles are externally tangent.

The purpose of the CC module is to reduce this contact error by allowing both centers and both radii to change simultaneously.
Let 
\[
\mathbf u=\frac{\mathrm{P}_2-\mathrm{P}_1}{\|\mathrm{P}_2-\mathrm{P}_1\|}
\]
be the unit vector pointing from $\mathrm{P}_1$ toward $\mathrm{P}_2$, and this vector is well defined.

Let $q$ be the contraction factor for the contact error. 
The correction is distributed equally among the two center positions and the two radii. 
Set \[
\alpha=\frac{1-q}{4}.
\]  
One application of the CC module is then defined by 
\begin{align*}
\mathrm{P}_1' &= \mathrm{P}_1+\alpha d \mathbf u, \qquad \mathrm{P}_2' = \mathrm{P}_2-\alpha d \mathbf u, \\
r_1' &= r_1+\alpha d, \qquad r_2'=r_2+\alpha d.
\end{align*}

Since the two centers each move by $\alpha d$ toward each other, the new center distance is
$R'=R-2\alpha d$.
At the same time, the sum of the new radii is $r_1'+r_2' = r_1+r_2+2\alpha d$.

Therefore the new contact error is 
\begin{align*}
d' =R'-r_1'-r_2'=(R-2\alpha d) -(r_1+r_2+2\alpha d) =d-4\alpha d =qd.
\end{align*}
Hence one application of the CC module multiplies the contact error by $q$.

In the same manner as in the LC module, no separate case distinction of  the sign of $d$ is necessary. 
If $d_0$ is the initial contact error, then after $n$ applications of the module, $d_n=q^nd_0$. Thus, for $|q|<1$, we have $d_n\to 0$ as $n\to\infty$.
The two circles converge asymptotically to an externally tangent configuration. 

\begin{figure}[hbt]
\begin{center}
\includegraphics[height=5cm]{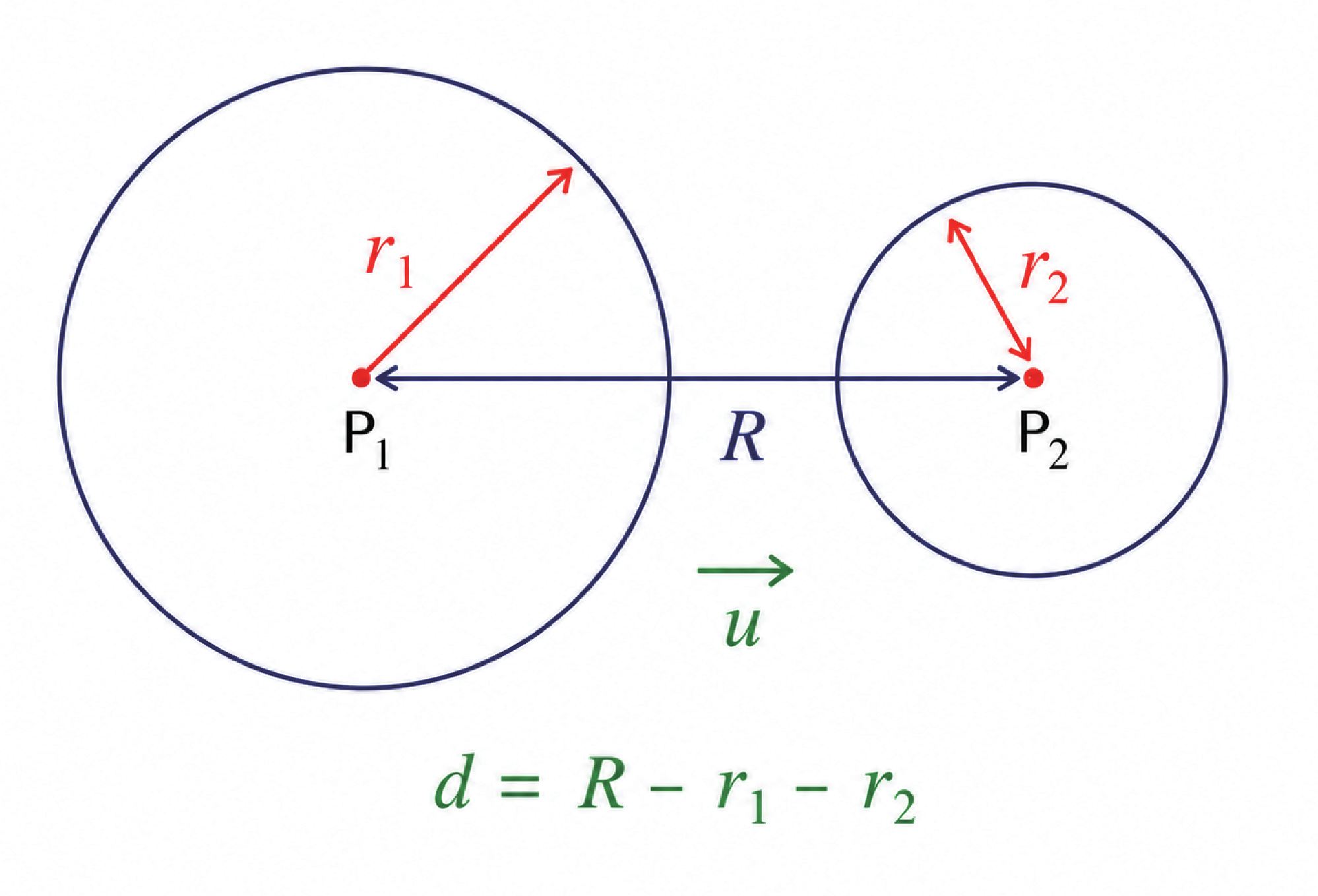}
\caption{CC Module}
\label{fig:CCmodule}
\end{center}
\end{figure}

\subsection{Markov process for the Malfatti circles}
We are now ready to define a {\bf Markov process for the Malfatti circles} using the LC and CC modules introduced above.

Let \(\mathrm{ABC}\) be a fixed triangle, and let \(C_1,C_2,C_3\) be three circles whose centers and radii are allowed to vary. A state of the system is specified by the positions of the three centers together with the three radii. To fix the notation, we assume that the \(i\)-th circle initially has center \(\mathrm{P}_i^{(0)}\) and radius \(r_i^{(0)}\). After applying the Markov process \(n\) times, we denote its center and radius by \(\mathrm{P}_i^{(n)}\) and \(r_i^{(n)}\), respectively. Thus, the state at step \(n\) may be represented by a vector
\[
X_n=(\mathrm{P}_1^{(n)}, r_1^{(n)}, \mathrm{P}_2^{(n)}, r_2^{(n)}, \mathrm{P}_3^{(n)}, r_3^{(n)}) \in \mathbb{R}^9. 
\]

For the Malfatti configuration $X_{\mathrm{Mal}}$, each circle must be tangent to two sides of the triangle, and each pair of circles must be externally tangent. Accordingly, we introduce six LC modules,
\[
(C_1,\mathrm{AB}),\quad (C_1,\mathrm{AC}),\quad
(C_2,\mathrm{AB}),\quad (C_2,\mathrm{BC}),\quad
(C_3,\mathrm{BC}),\quad (C_3,\mathrm{CA}),
\]
and three CC modules,
\[
(C_1,C_2),\qquad (C_2,C_3),\qquad (C_1,C_3).
\]

At each step, all nine modules act simultaneously on the current configuration. Each module computes its own contact error from the present state and produces a small correction of the relevant center positions and radii. The total displacement of each center and the total change of each radius are obtained by adding the corrections contributed by all modules involving that circle. This determines the next state \(X_{n+1}\) uniquely from the current state \(X_n\).

Thus, the evolution can be written in the folowing form
\[
X_{n+1}=T_q(X_n),
\]
where \(T_q\) is the update map determined by the six LC modules and the three CC modules, and \(q\) is the contraction parameter used in each module.

Since the next state depends only on the current state, not on the system's history, this iteration may be regarded as a discrete-time Markov process; more precisely, in the present deterministic setting, it is a deterministic Markov process.

A Malfatti configuration $X_{\mathrm{Mal}}$ is a fixed point of this process, that is, $T_q(X_{\mathrm{Mal}})=X_{\mathrm{Mal}}$. Indeed, when the three circles are the Malfatti circles, all six line-circle contact errors and all three circle-circle contact errors are zero. 

The main question is therefore whether repeated application of \(T_q\), starting from an arbitrary or suitably restricted initial configuration, converges to this fixed point.

\section{Proofs and discussions}

\subsection{Problem 1}

Our original first problem was simple to state: {\sl Is the Malfatti configuration a stable fixed point of the iteration?}

Once the six LC modules and three CC modules have been introduced, the Malfatti configuration is immediately seen to be a fixed point, since all line-circle and circle-circle contact errors vanish there. The more substantial issue is stability. 
Thus, the problem is naturally reformulated as a local stability problem for a discrete dynamical system. The aim is to determine whether the Malfatti configuration is a locally asymptotically stable fixed point of the update map $T_q$ defined by the LC and CC modules.

\subsubsection{Notation and Theorem}

The preceding discussion shows that the Malfatti configuration is not only a fixed point of the iteration $T_q$ defined by the six LC modules and the three CC modules, but is in fact locally asymptotically stable for $q=0.9$. 

\begin{theorem} (Local Stability of the Malfatti Configuration)

Let \(ABC\) be a nondegenerate triangle, and let \(X_{\mathrm{Mal}}\) denote its Malfatti configuration. Consider the discrete-time dynamical system defined by the simultaneous action of the nine modules as in 1.5 at each iteration.
Let \(T_{0.9}\) denote the resulting simultaneous update map. Then there exists a sufficiently small neighborhood \(U\) of \(X_{\mathrm{Mal}}\) such that, for every initial configuration $X_0\in U$, the sequence defined recursively by $X_{n+1}=T_{0.9}(X_n)$ converges geometrically to \(X_{\mathrm{Mal}}\). 
Hence \(X_{\mathrm{Mal}}\) is a locally asymptotically stable fixed point of the Malfatti iteration $T_{0.9}$.
\end{theorem}

\subsubsection{Proof by ChatGPT}

We now show that the simultaneous action of the six LC modules and the three CC modules can be interpreted as a gradient-descent step for a natural error functional.

Let $X=(\mathrm{P}_1,r_1,\mathrm{P}_2,r_2,\mathrm{P}_3,r_3)\in \mathbb R^9$ denote the state of the three-circle system. Here \(\mathrm{P}_i\in\mathbb R^2\) is the center of the \(i\)-th circle and \(r_i\in\mathbb R\) is its radius.  Let $l_1(X),\ldots,l_6(X)$ denote the six signed line-circle contact errors associated with the six LC modules, and let $c_{12}(X), c_{23}(X), c_{13}(X)$ denote the three signed circle-circle contact errors associated with the three CC modules. 
Each of these quantities vanishes precisely when the corresponding tangency condition is satisfied.

For the choice $q=0.9$, one application of the $k$-th LC module produces the correction
\begin{equation}
 \Delta X_k^{\mathrm{LC}} = -\frac{1}{20} l_k(X)\nabla l_k(X). \label{2.1.1}
\end{equation}
Here the gradient $\nabla l_k(X)$ is taken with respect to all nine coordinates of the state vector \(X\). Remark that $\frac{1}{20} = (1-0.9)/2.$ Since \(l_k\) depends only on the center and radius of the circle involved in that LC constraint, the remaining components of \(\nabla l_k\) are zero.

Similarly, the CC module corresponding to the pair \(C_i,C_j\) produces the correction
\[ \Delta X_{ij}^{\mathrm{CC}} = -\frac{1}{40}c_{ij}(X)\nabla c_{ij}(X)\]
for $(i,j) \in \{ (1,2), (2,3), (1,3)\}.$

Remark that $\frac{1}{40} = (1-0.9)/4.$ Again, the gradient $\nabla c_{ij}(X)$ is regarded as a vector in the full nine-dimensional state space, with zero components in the variables not involved in the corresponding CC constraint.

The important point is that all nine corrections are evaluated from the same current state \(X\). In other words, the modules are applied simultaneously rather than sequentially: no module uses a state that has already been modified by another module during the same iteration. Consequently, the total correction is obtained simply by adding the nine individual corrections. Thus,
\[
T_q(X)-X
=
-\frac{1}{20}\sum_{k=1}^{6}
l_k(X)\nabla l_k(X) - \frac{1}{40}\sum_{(i,j)}
c_{ij}(X)\nabla c_{ij}(X)
\]
Now define the error functional
\[
E_q(X) = \frac{1}{40}\sum_{k}l_k(X)^2
+ \frac{1}{80}
\sum_{(i,j)} c_{ij}(X)^2
\]
By the chain rule,
it follows that the simultaneous update map can be written in the compact form
\[ T_q(X)=X-\nabla E_q(X). \]

Thus, the Malfatti iteration is precisely a gradient-descent iteration, with unit step size, for the weighted sum of squared contact errors \(E\).

Since \(E_q(X)\ge 0\), and since $E_q(X)=0$ if and only if all six LC contact errors and all three CC contact errors vanish simultaneously, the zero set of \(E_q\) consists exactly of configurations satisfying all nine prescribed tangency conditions. In particular, the known Malfatti configuration is a global minimizer of \(E_q\), with minimum value \(0\).


Differentiating \(T_q\), we obtain the Hessian at the Malfatti configuration,
\[
\nabla^2T_q(X_{\mathrm{Mal}})=I-\nabla^2E_q(X_{\mathrm{Mal}}).
\]
We now compute the Hessian of \(E_q\) at \(X_{\mathrm{Mal}}\). 
For a typical LC term $\frac{1}{40}l_k(X)^2$, its Hessian is
\[
\frac{1}{20}
\left(
\nabla l_k\,\nabla l_k^{T} + l_k\,\nabla^2l_k
\right).
\]
Since \(l_k(X_{\mathrm{Mal}})=0\), the second term vanishes at \(X_{\mathrm{Mal}}\), and hence
\[
\nabla^2
\left(
\frac{1}{40}l_k^2
\right)(X_{\mathrm{Mal}})
=\frac{1}{20}
\nabla l_k(X_{\mathrm{Mal}})\nabla l_k(X_{\mathrm{Mal}})^{T}.
\]

Similarly, for a CC term $\frac{1}{80}c_{ij}(X)^2$, we have
\[
\nabla^2
\left(
\frac{1}{80}c_{ij}^2
\right)(X_{\mathrm{Mal}})
=
\frac{1}{40}
\nabla c_{ij}(X_{\mathrm{Mal}})\nabla c_{ij}(X_{\mathrm{Mal}})^{T}.
\]

Consequently, $\nabla^2T_q(X_{\mathrm{Mal}})=I-H$, where
\begin{equation}
H=
\frac{1}{20}
\sum_{\mathrm{LC}}
\nabla l_k(X_{\mathrm{Mal}})\nabla l_k(X_{\mathrm{Mal}})^{T}
+
\frac{1}{40}
\sum_{\mathrm{CC}}
\nabla c_{ij}(X_{\mathrm{Mal}})\nabla c_{ij}(X_{\mathrm{Mal}})^{T}.
\label{2.1.2}
\end{equation}

The matrix \(H\) is symmetric. Moreover, it is positive semidefinite. Indeed, for any vector \(v\in\mathbb R^9\),
\[
v^{T}Hv
=
\frac{1}{20}
\sum_{\mathrm{LC}}
\left(
\nabla l_k(X_{\mathrm{Mal}})^{T}v
\right)^2
+
\frac{1}{40}
\sum_{\mathrm{CC}}
\left(
\nabla c_{ij}(X_{\mathrm{Mal}})^{T}v
\right)^2
\ge 0.
\]

\bigskip

The crucial issue is whether \(H\) is in fact positive definite. That is, the linear independence of the nine contact constraints at the Malfatti configuration is the central step in establishing the local stability of \(X_{\mathrm{Mal}}\).

Let the angles of the triangle be denoted by $A, B, C$, and let the side lengths be $c=AB, a=BC, b=CA$.
Assume that all six LC conditions are satisfied.  Then the center of \(C_1\) lies on the internal angle bisector of \(A\), the center of \(C_2\) lies on the internal angle bisector of \(B\), and the center of \(C_3\) lies on the internal angle bisector of \(C\).  Therefore, once the six LC conditions are imposed, the remaining degrees of freedom are essentially the three radii $r_1, r_2, r_3$. 

Now \(C_1\) and \(C_2\) are both tangent to the side \(AB\). Suppose that two circles of radii \(r_1\) and \(r_2\) are tangent to the same line and are externally tangent to each other. Then the distance, measured along the line, between the orthogonal projections of their centers onto the line is $2\sqrt{r_1r_2}$.  Hence,
\begin{equation}
c
=
r_1\cot\frac{A}{2}
+
2\sqrt{r_1r_2}
+
r_2\cot\frac{B}{2}.
\label{2.1.3}
\end{equation}

Similarly, we obtain
\begin{equation}
a
=
r_2\cot\frac{B}{2}
+
2\sqrt{r_2r_3}
+
r_3\cot\frac{C}{2},
\label{2.1.4}
\end{equation}
and
\begin{equation}
b
=
r_3\cot\frac{C}{2}
+
2\sqrt{r_3r_1}
+
r_1\cot\frac{A}{2}.
\label{2.1.5}
\end{equation}

These three equations are precisely the three CC conditions that remain after the six LC conditions have been imposed.

For simplicity, set
\(
\alpha=\cot\frac{A}{2}, 
\beta=\cot\frac{B}{2}, 
\gamma=\cot\frac{C}{2}.
\)
Since \(A,B,C\in(0,\pi)\), we have
\(
\alpha>0, \beta>0, \gamma>0.
\)
Let \(F_1,F_2,F_3\) denote the radius-dependent expressions appearing in
the right-hand sides of equations (\ref{2.1.3})--(\ref{2.1.5}), respectively. 
The Jacobian matrix of \((F_1,F_2,F_3)\) with respect to
\((r_1,r_2,r_3)\) is
\[
\frac{\partial(F_1,F_2,F_3)}
     {\partial(r_1,r_2,r_3)}
=
\begin{pmatrix}
\displaystyle
\alpha+\sqrt{\frac{r_2}{r_1}}
&
\displaystyle
\beta+\sqrt{\frac{r_1}{r_2}}
&
0
\\[3mm]
0
&
\displaystyle
\beta+\sqrt{\frac{r_3}{r_2}}
&
\displaystyle
\gamma+\sqrt{\frac{r_2}{r_3}}
\\[3mm]
\displaystyle
\alpha+\sqrt{\frac{r_3}{r_1}}
&
0
&
\displaystyle
\gamma+\sqrt{\frac{r_1}{r_3}}
\end{pmatrix}.
\]

For convenience, write this matrix in the form
\[
\begin{pmatrix}
u_{11}&u_{12}&0\\
0&u_{22}&u_{23}\\
u_{31}&0&u_{33}
\end{pmatrix}.
\]
Since \(r_1,r_2,r_3>0\) and \(\alpha,\beta,\gamma>0\), all six entries are strictly positive, that is,
\[u_{11}, u_{12}, u_{22}, u_{23}, u_{31}, u_{33} > 0.\]

Its determinant is $u_{11}u_{22}u_{33}+ u_{12}u_{23}u_{31}>0$, hence the Jacobian matrix is nonsingular. 
The six LC conditions reduce the original nine-dimensional state space to three remaining
degrees of freedom, which may be parametrized locally by $r_1, r_2, r_3$.
The nonsingularity of the above Jacobian shows that the three CC
conditions eliminate these remaining three infinitesimal degrees of
freedom.

Consequently, the gradients of the six LC constraints together with
those of the three CC constraints span the full state space $\mathbb{R}^9$.
Equivalently, the nine constraint gradients are linearly independent at
the Malfatti configuration, and \(H\) is positive definite.

\bigskip

We next need to verify that a single gradient-descent step is not so large as
to cause divergence.

For an LC contact error $l=r-h$, where \(h\) denotes the distance from the center of the circle to the
corresponding side of the triangle, we have $\|\nabla l\|^2=2$.
Indeed, the gradient is $(-u_x, -u_y, 1)$, where $(u_x, u_y)$ is the internal unit normal vector of the side.

Recall that the contribution of one LC constraint to \(H\) is $\frac{1}{20}\nabla l\,\nabla l^{T}$.
Since a rank-one matrix of the form \(vv^{T}\) has largest eigenvalue
\(\|v\|^2\), the largest eigenvalue of this contribution is $\frac{1}{20}\cdot 2=0.1$.
There are six LC constraints. Hence, using the subadditivity of the largest
eigenvalue for symmetric positive semidefinite matrices, the total
contribution of the six LC terms has largest eigenvalue at most $0.6$.

For a CC contact error, $c_{ij} = \|\mathrm{P}_i-\mathrm{P}_j\|-r_i-r_j$, the gradient has four nonzero components: all of the derivatives with respect to \(\mathrm{P}_i\), \(\mathrm{P}_j\), \(r_i\) and \(r_j\) both have absolute value \(1\). Therefore,
$\|\nabla c_{ij}\|^2=4$. The contribution of one CC constraint to \(H\) is
$\frac{1}{40}\nabla c_{ij}\,\nabla c_{ij}^{T}$, and hence its largest eigenvalue is $0.1$.
Since there are three CC constraints, their total contribution has largest
eigenvalue at most $0.3$.

Combining the LC and CC contributions, we obtain the largest eigenvalue $\lambda_{\max}(H)$ of $H$ satisfy 
\[ \lambda_{\max}(H) \le 0.6+0.3=0.9. \] On the other hand, as proved above, \(H\) is positive definite. Hence its
smallest eigenvalue satisfies $\lambda_{\min}(H)>0$.

Since $\nabla^2 T_q(X_{\mathrm{Mal}})=I-H$, the any eigenvalues $\lambda$ of \(\nabla^2 T_q(X_{\mathrm{Mal}})\) satisfies 
\[ 0.1 \leq \lambda < 1\]
Consequently, the spectal radius is less than $1$, and $\|\nabla^2 T_q(X_{\mathrm{Mal}})\|_2 < 1$.

\bigskip

By continuity of \(\nabla^2 T_q\), there therefore exist a sufficiently small
neighborhood \(U\) of \(X_{\mathrm{Mal}}\) and a constant \(q_0\) with $ 0<q_0<1$ such that
\[
\|\nabla^2T_q(X)\|_2\leq q_0
\qquad
\text{for all } X\in U.
\]
Shrinking \(U\), if necessary, we may assume that it is convex.

For any \(X_0\in U\), the mean value estimate then gives
\[
\|T_q(X_0)-T_q(X_{\mathrm{Mal}})\| \leq q_0\|X_0-X_{\mathrm{Mal}}\|.
\]
Since \(T_q(X_{\mathrm{Mal}})=X_{\mathrm{Mal}}\), this becomes
\[
\|T_q(X_0)-X_{\mathrm{Mal}}\| \leq q_0\|X-X_{\mathrm{Mal}}\|.
\]
All iterates remain in this neighborhood and satisfy
\[
\|T_q^n(X_0)-X_{\mathrm{Mal}}\| \leq q_0^n\|X_0-X_{\mathrm{Mal}}\|\longrightarrow 0.
\]
for $n\to\infty$. Hence
\[
T_q^n(X_0)\longrightarrow X_{\mathrm{Mal}}.
\]
Thus the Malfatti configuration \(X_{\mathrm{Mal}}\) is a locally asymptotically stable
fixed point of the simultaneous LC--CC iteration. Moreover, the convergence
is geometric.

\subsection{Problem 2}

We now turn to the second problem, which concerns the range of admissible
values of the parameter \(q\).

For a single LC or CC module, the corresponding contact error \(d\) is
updated according to
$
d' = qd.
$
Hence, if the module is considered in isolation, the contact error converges
to zero under repeated application whenever
$
|q|<1.
$

It is therefore natural to ask whether the simultaneous Malfatti iteration,
consisting of six LC modules and three CC modules, remains stable throughout
this entire range. Accordingly, in Problem 2 we investigate which values of \(q\) actually guarantee local convergence to the Malfatti configuration.

\subsubsection{Notation and Theorem}

Let $X_{\mathrm{Mal}}$ denote the Malfatti configuration of a fixed nondegenerate
triangle $ABC$. At $X_{\mathrm{Mal}}$, define the symmetric matrix
\[
 M
 =
 \frac{1}{2}\sum_{\mathrm{LC}}
 \nabla l_k(X_{\mathrm{Mal}})\nabla l_k(X_{\mathrm{Mal}})^{T}
 +
 \frac{1}{4}\sum_{\mathrm{CC}}
 \nabla c_{ij}(X_{\mathrm{Mal}})\nabla c_{ij}(X_{\mathrm{Mal}})^{T},
\]
where $l_k$ are the six line--circle contact errors and $c_{ij}$ are the
three circle--circle contact errors. As shown above, the nine contact
constraints are locally independent at the Malfatti configuration;
hence $M$ is positive definite.

The following theorem gives both a triangle-dependent stability
criterion and a uniform criterion independent of the shape of the
triangle.

\begin{theorem}[Stability range of the contraction parameter]
Let $ABC$ be a nondegenerate triangle, let $X_{\mathrm{Mal}}$ be its Malfatti
configuration, and let $T_q$ be the simultaneous update map obtained
from the six LC modules and the three CC modules, where each individual
module multiplies its contact error by a factor $q$, with $-1<q<1$.

(1) Let $\lambda_{\max}(M)$ denote the largest eigenvalue of the matrix $M$
defined above. Then the Malfatti configuration $X_{\mathrm{Mal}}$ is a locally
asymptotically stable fixed point of $T_q$ whenever
\[
1-\frac{2}{\lambda_{\max}(M)}<q<1.
\]

(2) Furthermore, one has the triangle-independent uniform condition
\[
\frac{1}{2}<q<1
\]
guarantees local asymptotic stability of the Malfatti configuration for
every nondegenerate triangle.

\end{theorem}

\subsubsection{Proof by ChatGPT}

(1) 
We replace the value $0.9$ used above by a general parameter $q$.
Thus, for each individual LC or CC module, the corresponding contact
error $d$ is transformed according to
\[
    d \longmapsto qd,
    \qquad -1<q<1.
\]
For an isolated module, this condition is sufficient to guarantee
convergence of the contact error to zero.

Let $X\in\mathbb R^9$ denote the state of the three-circle system, consisting of the
three center positions and the three radii. For an LC constraint with
contact error $l$, the update is
\[
    X \longmapsto
    X-\frac{1-q}{2}\,l\,\nabla l.
\]
Similarly, for a CC constraint with contact error $c$, the update is
\[
    X \longmapsto
    X-\frac{1-q}{4}\,c\,\nabla c.
\]

Since the six LC modules and the three CC modules are applied
simultaneously, their corrections are added. Hence the complete
update map can be written as
\[
    T_q(X)=X-(1-q)\nabla F(X),
\]
where
\[
    F(X)
    =
    \frac{1}{4}\sum_{\mathrm{LC}} l_k(X)^2
    +
    \frac{1}{8}\sum_{\mathrm{CC}} c_{ij}(X)^2.
\]
Thus, the Malfatti process may be regarded as a gradient descent method
for the contact-error energy $F$, with step size $1-q$.

Let $X_{\mathrm{Mal}}$ be the Malfatti configuration. Since all contact errors vanish
at $X_{\mathrm{Mal}}$,
\[
    l_k(X_{\mathrm{Mal}})=0,
    \qquad
    c_{ij}(X_{\mathrm{Mal}})=0,
\]
the linearization of $T_q$ at $X_{\mathrm{Mal}}$ is
\[
    \nabla^2 T_q(X_{\mathrm{Mal}})
    =
    I-(1-q)M,
\]
where
\[
    M
    =
    \frac{1}{2}
    \sum_{\mathrm{LC}}
    \nabla l_k(X_{\mathrm{Mal}}) \nabla l_k(X_{\mathrm{Mal}})^{T}
    +
    \frac{1}{4}
    \sum_{\mathrm{CC}}
    \nabla c_{ij}(X_{\mathrm{Mal}}) \nabla c_{ij}(X_{\mathrm{Mal}})^{T}.
\]

As discussed above, the nine contact constraints are locally independent
at the Malfatti configuration. Consequently, the matrix $M$ is positive
definite. Let its eigenvalues be
\[
    0<\lambda_1\leq\lambda_2\leq\cdots\leq\lambda_9
    =\lambda_{\max}(M).
\]
The eigenvalues of the linearized update map $\nabla^2 T_q(X_{\mathrm{Mal}})$ are therefore
\[
    1-(1-q)\lambda_i,
    \qquad i=1,\ldots,9.
\]

Local asymptotic stability requires that all these eigenvalues have
absolute value strictly less than one:
\[
    \left|1-(1-q)\lambda_i\right|<1
    \qquad
    (i=1,\ldots,9).
\]
Since $q<1$ and all $\lambda_i$ are positive, this condition is equivalent
to
\[
    (1-q)\lambda_{\max}(M)<2.
\]
Hence the exact local stability condition for a fixed triangle is
\[
    q>    1-\frac{2}{\lambda_{\max}(M)}.
\]

Thus, once the geometry of the triangle is fixed, the lower stability
threshold can in principle be determined by computing the largest
eigenvalue of $M$.

For example, if ABC is the regular triangle, $\lambda_{\max}=\frac{3+\sqrt{5-2\sqrt{3}}}{2} \sim 2.11965684$.
Thus the critical value of $q$ is $\frac{-1+\sqrt{5-2\sqrt{3}}}{3+\sqrt{5-2\sqrt{3}}} \sim 0.056451$.

\bigskip

(2) We next derive a uniform sufficient condition that does not depend on
the shape of the triangle. Let
\[
    v=(\mathbf p_1,s_1,\mathbf p_2,s_2,\mathbf p_3,s_3)
\]
be an infinitesimal variation of the state, where \(\mathbf p_i\in\mathbb{R}^2 \) represents the displacement of the center of the $i$-th circle and
$s_i\in\mathbb{R}$ represents the corresponding variation of its radius.

Consider first a single LC term. If $\mathbf n$ is the unit normal vector to
the corresponding side of the triangle, its contribution to the
quadratic form $\mathbf v^TM\mathbf v$ is \(\frac{1}{2}  \left(s_i-(\mathbf n, \mathbf p_i)\right)^2 \).
Using $\frac{1}{2}(a-b)^2\leq a^2+b^2$, we obtain
\[
    \frac{1}{2}
    \left(s_i - (\mathbf n, \mathbf p_i)\right)^2
    \leq
    s_i^2+(\mathbf n, \mathbf p_i)^2
    \leq
    s_i^2+\| \mathbf p_i\|^2.
\]
Since each circle is involved in exactly two LC constraints, the total
contribution of the six LC terms satisfies
\[
    \text{LC contribution}
    \leq
    2\sum_{i=1}^{3}
    \left(
        \| \mathbf p_i\|^2+s_i^2
    \right).
\]

Now consider a CC constraint between circles $i$ and $j$.
Let $u$ be the unit vector joining their centers at the Malfatti
configuration. The corresponding contribution to $\mathbf v^TM\mathbf v$ is
\[
    \frac{1}{4}
    \left\{
        (\mathbf u,(\mathbf p_i-\mathbf p_j))-s_i-s_j
    \right\}^2.
\]
Applying the elementary inequality $(a+b+c+d)^2 \leq 4(a^2+b^2+c^2+d^2)$, we obtain
\[
\begin{split}
    \frac{1}{4}
    \left\{
        (\mathbf u,(\mathbf p_i-\mathbf p_j))-s_i-s_j
    \right\}^2
    &\leq
    (\mathbf u, \mathbf p_i)^2
    +(\mathbf u, \mathbf p_j)^2
    +s_i^2+s_j^2
    \\
    &\leq
    \| \mathbf p_i\|^2
    +\| \mathbf p_j\|^2
    +s_i^2+s_j^2.
\end{split}
\]
There are three CC constraints, and each circle occurs in exactly two
of them. Therefore,
\[
    \text{CC contribution}
    \leq
    2\sum_{i=1}^{3}
    \left(
        \| \mathbf p_i\|^2+s_i^2
    \right).
\]
Combining the LC and CC estimates gives
\[
    \mathbf v^TM\mathbf v
    \leq
    4\sum_{i=1}^{3}
    \left(
        \|\mathbf p_i\|^2+s_i^2
    \right)
    =
    4\| \mathbf v\|^2.
\]
It follows that
\[
    \lambda_{\max}(M)\leq 4.
\]

Substituting this estimate into the triangle-dependent stability
condition yields
\[
    1-\frac{2}{\lambda_{\max}(M)}
    \leq
    1-\frac{2}{4}
    =
    \frac{1}{2}.
\]
Consequently, $ \frac{1}{2}<q<1 $ is a sufficient condition for local asymptotic stability of the
Malfatti configuration for every nondegenerate triangle.

\section{Discussion with ChatGPT}

\subsection{Related research on algrithm}
\label{sec:related}

The author asked ChatGPT to investigate the prior research related to this study. Specifically, rather than asking about research motivation, he asked whether any papers presented proof frameworks or similar problem formulations and calculations.

The first paper is \cite{Witkin} by Witkin et. al.  This paper describes a method that represents geometric constraints in CG space as non-negative energy functionals, sums them, and moves parameters in the direction of the energy gradient to converge on a configuration that satisfies the constraints. Furthermore, the variables are not limited to mere point coordinates; the system can adjust any model parameter. 

The second is \cite{Muller} by M\"uller et. al.  In computer graphics—particularly in real-time physics simulation—this method is used to animate deformable objects such as cloth stably and rapidly, as well as to handle collisions between objects. Technically speaking, if the physical constraints are expressed as $C(p) = 0$, they made a local correction by the formula
\[\Delta p=-\frac{C(p)}{\| \nabla C(p)\|^2 }\nabla C(p).
\]
This correction method bears a striking resemblance to the equation (\ref{2.1.1})  in this paper.  ChatGPT said  ''Rather than viewing this as a coincidence, it is better to consider that we naturally arrived at the same equation based on the requirement to distribute the constraint error symmetrically across the relevant degrees of freedom. The PBD paper also presents this equation as a Newton–Raphson step for a single constraint.''

\subsection{Related research on circle packing}

Collins et al. [4] proposed an iterative algorithm for constructing circle packings with a prescribed tangency pattern by repeatedly adjusting the radii of the circles. They presented efficient numerical algorithms for approximating such circle packings in both Euclidean and hyperbolic geometries.
The underlying circle-packing problem is based on the classical circle packing theorem, originating with Koebe and later brought to renewed attention by Thurston around 1977. The theorem states that every finite planar graph can be realized as the tangency graph of a collection of circles.
Although the algorithm of Collins et al. is also based on iterative adjustment toward prescribed tangency conditions, its mechanism differs from that of the CC module introduced in the present paper. In particular, our CC module directly modifies both the centers and the radii of a pair of circles according to their local tangency error.

\subsection{More Malfatti-like circles}

\begin{figure}[hbt]
\begin{center}
\includegraphics[height=5cm]{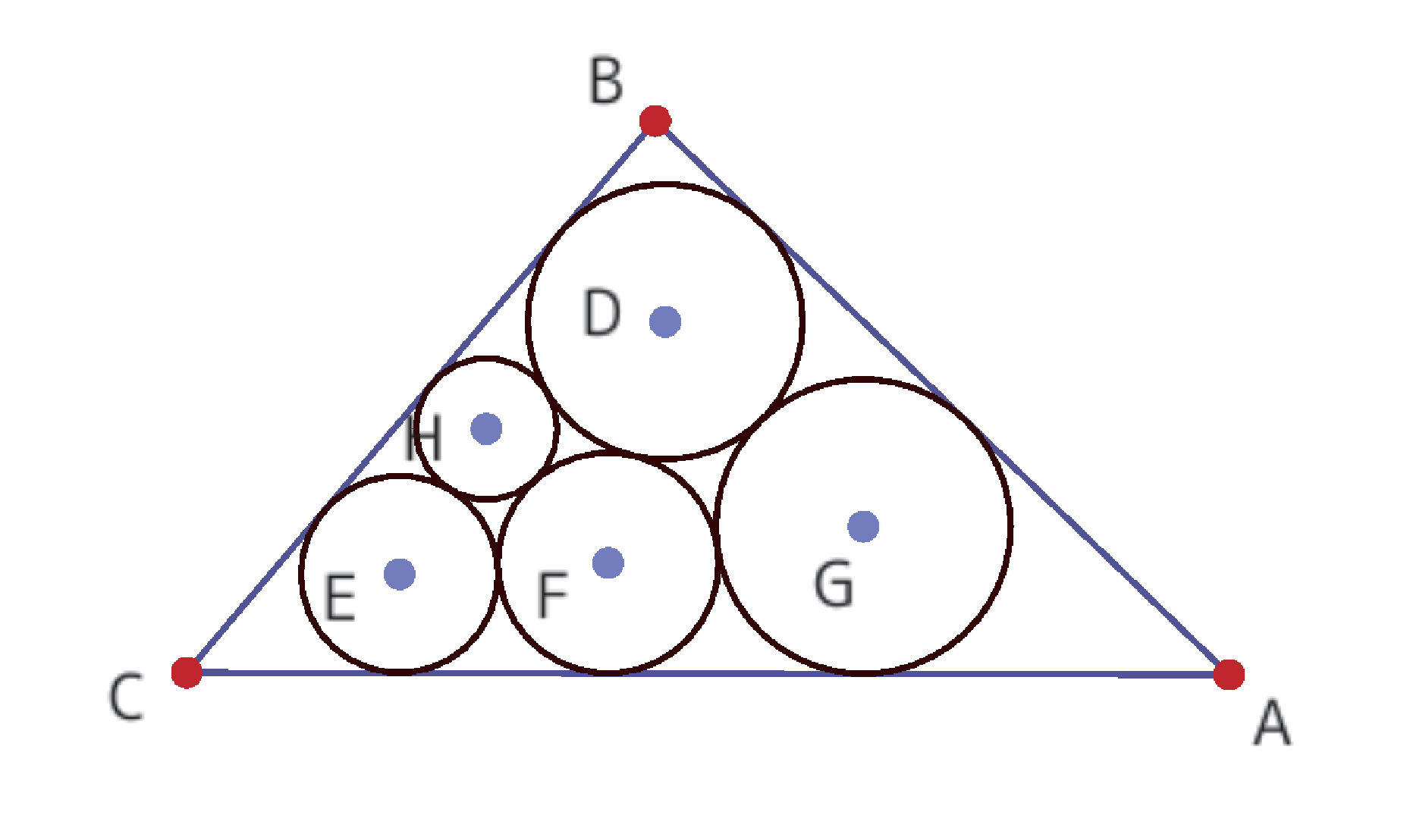}
\caption{Malfatti-like circles}
\label{fig:Malfatti-like-circles}
\end{center}
\end{figure}

We can develop a similar stability analysis for more complex circular configurations, such as the one shown in Figure \ref{fig:Malfatti-like-circles}. For each prescribed tangency between a circle and a side of the triangle, we introduce an LC module, and for each prescribed tangency between two circles, we introduce a CC module. All modules are then applied simultaneously.

If \(N\) circles are involved, the state space has dimension \(3N\), since each circle has two center coordinates and one radius. Let \(X^*\) denote a configuration satisfying all prescribed tangency conditions. As in the Malfatti case, we can write the simultaneous update in the form.
\[
T_q(X)=X-(1-q)\nabla F(X),
\]
where \(F\) is the weighted sum of the squared LC and CC contact errors.

At the fixed point \(X^*\), the linearization is
\[
\nabla^2T_q(X^*)=I-(1-q)M,
\]
where \(M\) is a positive semidefinite matrix determined by the gradients of all tangency constraints. If these constraints are locally independent so that
\(
M>0,
\)
then \(X^*\) is an isolated fixed point, and the same argument as for the Malfatti configuration gives the local stability condition
\[
\left|1-(1-q)\lambda_i(M)\right|<1
\]
for every eigenvalue of \(M\). Equivalently,
\[
q>1-\frac{2}{\lambda_{\max}(M)}.
\]
Thus, the general method extends almost immediately to configurations with more circles. The main additional task for each particular pattern, such as the one in Figure \ref{fig:Malfatti-like-circles}, is to verify that the prescribed LC and CC constraints are locally independent and to estimate or compute \(\lambda_{\max}(M)\). The resulting lower bound for \(q\) will generally depend on the combinatorial pattern of tangencies as well as on the geometry of the particular configuration.

\subsection{Problem 3(initial value problem)}
A global convergence result for arbitrary initial configurations appears to be much more difficult. In particular, the simultaneous interaction of many LC and CC modules may produce complicated trajectories, and there is no obvious invariant region or global Lyapunov argument that guarantees convergence to the desired configuration. Therefore, at present, a satisfactory solution to the global initial-value problem seems out of reach.

\subsubsection{Case study for all circles are the triangle and have pairwise disjoint}
One may ask whether a global initial-value result can be obtained under a natural geometric assumption, for example, that all three circles are initially contained in the triangle and have pairwise disjoint interiors. However, this condition is not preserved by the simultaneous LC–CC update. In fact, even for initial configurations arbitrarily close to the Malfatti configuration, it is possible to choose three initially non-overlapping circles such that, after a single iteration, some pair of circles overlaps; equivalently, although initially $c_{ij}>0$, one may have $c_{ij}'<0$ after one step. Thus, the set defined by non-overlapping circles is not an invariant set of the process, making it difficult to derive a global convergence theorem based solely on this initial condition.

\subsubsection{Non-stable critical configuration}
\begin{figure}[hbt]
\begin{center}
\includegraphics[height=5cm]{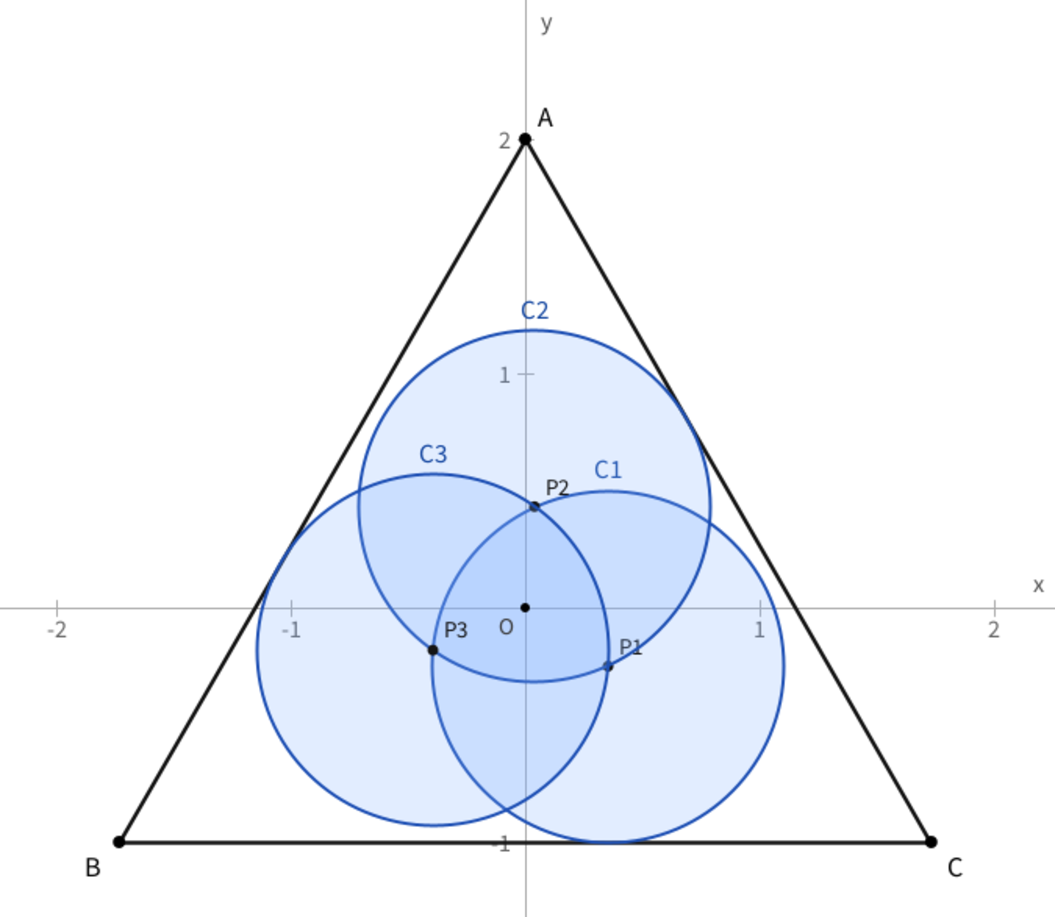}
\caption{non-stable critical configuration}
\label{fig:non-stable}
\end{center}
\end{figure}

Besides the Malfatti configuration, the energy functional \(H\), as in formula (\ref{2.1.2}), may have other critical points. For example, consider the equilateral triangle
\[
A=(0,2),\qquad B=(-\sqrt{3},-1),\qquad C=(\sqrt{3},-1).
\]
A direct calculation shows that \(H\) has a non-Malfatti critical point in which all three circles have radius
\[
r_1=r_2=r_3=\frac34,
\]
and their centers are
\begin{align*}
P_1 &=
\left(
\frac{\sqrt2}{4},-\frac14
\right),\\
P_2 &=
\left(
\frac{\sqrt3-\sqrt2}{8},
\frac{1+\sqrt6}{8}
\right),\\
P_3 &=
\left(
-\frac{\sqrt3+\sqrt2}{8},
\frac{1-\sqrt6}{8}
\right).
\end{align*}
At this configuration, see Figure \ref{fig:non-stable}, $\nabla H=0$, although none of the prescribed tangency conditions is simultaneously satisfied. In particular,
\[
c_{12}=c_{23}=c_{13}=-\frac34,
\]
so the three circles overlap substantially. This provides an explicit example of a non-Malfatti critical point of \(H\); the corresponding critical point is unstable. There is also a mirror-image critical point obtained by reflection of the above configuration.

\subsubsection{Some more initial-value problem}
A natural question is whether the conditions 
\[ l_k>0 \quad (k=1,\ldots,6), \text{ and }c_{ij}>0 \quad (1\le i<j\le3)
\]
 are sufficient to guarantee convergence to the Malfatti configuration. Under our sign convention \(l_k=r-h\), the condition \(l_k>0\) means that the radius of the corresponding circle is larger than the distance from its center to the relevant side of the triangle, while \(c_{ij}>0\) means that the two circles \(C_i\) and \(C_j\) are separated and have disjoint interiors.

At present, we do not have a mathematical proof that these conditions imply convergence. However, in a numerical experiment with \(10{,}000\) randomly generated initial configurations satisfying these inequalities, every trajectory converged to the Malfatti configuration. Although such numerical evidence does not constitute a proof, it strongly suggests the following conjecture:
\[
l_k>0,\quad c_{ij}>0
\quad\Longrightarrow\quad
X_n\longrightarrow X_{\mathrm Mal}
\quad (n\to\infty).
\]
Thus, these inequalities appear to provide a promising candidate for a relatively large basin of attraction of the Malfatti configuration.

\end{document}